\documentclass[conference,a4paper]{IEEEtran}

\usepackage[T1]{fontenc}
\usepackage[utf8]{inputenc}

\usepackage{cite}
\usepackage[cmex10]{amsmath}
\usepackage{amssymb,amsfonts}
\usepackage{array}
\usepackage[caption=false,font=footnotesize]{subfig}
\usepackage{dblfloatfix}

\usepackage{algorithmic}
\usepackage{graphicx}
\graphicspath{{graphics/}{moregraphics/}}
\DeclareGraphicsExtensions{.pdf,.PDF,.jpeg,.JPEG,.jpg,.JPEG,.png,.PNG}

\usepackage{booktabs}

\usepackage{siunitx}

\usepackage{textcomp}
\usepackage{xcolor}
\usepackage[spaces,hyphens]{url}

\def\BibTeX{{\rm B\kern-.05em{\sc i\kern-.025em b}\kern-.08em
    T\kern-.1667em\lower.7ex\hbox{E}\kern-.125emX}}

\begin{document}

\title{Training without Gradients - A Filtering Approach}

\author{\IEEEauthorblockN{Isaac Yaesh}
\IEEEauthorblockA{\textit{Control Dept.} \\
\textit{Elbit Systems}\\
Ramat Hasharon, Israel \\
itzhak.yaesh@elbitsystems.com}
\and
\IEEEauthorblockN{Natan Grinfeld}
\IEEEauthorblockA{\textit{Control Dept.} \\
\textit{Elbit Systems}\\
Ramat Hasharon, Israel \\
nati.grinfeld@elbitsystems.com}

}
\maketitle

\begin{abstract}
  A particle filtering approach is suggested for the training of multi-layer neural networks without utilizing gradients calculation. The network weights are considered to be the components of the estimated state-vector of a noise driven linear system, whereas the neural network serves as the measurement function in the estimation problem. A simple example is used to provide a preliminary demonstration of the concept, which remains to be further studied for training deep neural networks.    
\end{abstract}

\section{Introduction}
Artificial Neural Networks (ANN) are motivated by the
nonlinear, complex and parallel structure of the human brain
\cite{Haykin:1998:NNC:521706} and are used in various applications such as control of
nonlinear systems, system identification and classification.
Traditional neural networks architectures include single and
multi layer feed-forward networks as well as recurrent networks
and they are based on interconnections of neurons of
the same type. The feed-forward networks are usually trained
by the back propagation algorithm (see e.g. \cite{Haykin:1998:NNC:521706}). However, back-propagation based training may face vanishing gradients as the result of the products of e.g. partial derivative of hyperbolic tangent functions which possess gradients in the range of $(0,1)$. A few solutions have been proposed to alleviate the vanishing gradient problem, such as using reLU activation functions, as well as pre-training one network level at a time and fine tuning using back propagation. A widely used approach to deal with the vanishing gradient problem is the one of ResNet which involves splitting the network to chunks and direct passes which actually make deep networks. ResNet indeed reduces much of the vanishing gradient problem but as they are in a way a cascade of shallow networks, they do not avoid the calculation of gradients. One approach which avoids calculation of gradients is the ADMM (Alternating Direction Method  of Multipliers \cite{admm}), where the network weights optimization problem is broken into a sequence of minimization sub-steps where each step can be globally optimized, however, originally not guaranteeing global convergence of the overall scheme. More recently such convergence was established \cite{admm1} under some conditions. Yet we would like to suggest the entirely different approach of Particle Filtering with implicit guarantee of global convergence and no need in calculating gradients. To this end, it is useful to recall a relation between filtering and back-propagation. The latter relation has been established in e.g. \cite{inproceedings_Hassibi}. It was shown there, that $H_\infty$-optimal filtering HF is closely related to back-propagation based training of neural networks. A similar result appears in \cite{hassibi1995h}.  Note that in HF the maximum ratio between the energy of the network output prediction errors and the energy of disturbances
(i.e. observation noise) is minimized as oppose to the Extended Kalman filter (EKF) which is a finite dimensional approximation to the optimal Minimal Mean Square Error (MMSE) estimate for the network weights. However, in the case of general multi-layer networks, gradients are still needed both in EKF and HF as it will be shown in the next section. In contrast,in Affine Functional Networks (AFN), a special case of Functional Neural Networks \cite{castillo1998introduction} with an affine relation between network output and its weights, the use of gradients is not needed. For this case, in \cite{yaesh2014robust} an $H_{\infty}$-filter has been derived, which was shown to be superior over the Kalman filter, being effective and a possible solution to the training of AFN. In AFN \cite{yaesh2014robust} the network weights are regarded as the components of the state-vector, which is estimated using either a standard Kalman filter or $H_\infty$ filter, both avoiding the use of gradients. Unfortunately, in most of the applications, multi-layer networks with a nonlinear relation rather than a linear or affine one between the network weights and output are used, and both the Kalman filter and $H_\infty$-optimal filters are to be replaced by approximations (e.g. Extended Kalman Filter - EKF) requiring the calculation of gradients.       

It is well known the MMSE is just a conditional expectation $\hat{x} =E\{x(t)/Y(t)\}$ of the state-vector $x(t)$ (i.e. the weights vector) on the history of the measurements $Y=\{y(\tau), \tau \leq t$ (i.e. training example outputs). In the general nonlinear networks case, the conditional expectation requires an infinite-dimensional filter \cite{jazwinski2007stochastic} which can be approximated in various ways, such as the EKF which can, however, even diverge in some special nonlinear cases. The most generally effective and widely used approximation of the MMSE is the particle filter (see e.g. \cite{Arulampalam02atutorial})
 (PF), where a finite-dimensional version of the conditional probability distribution function (PDF) is estimated where it is represented by points on the multidimensional graph known as particles.   

\section{Problem Formulation}
Consider the following training problem :
\begin{equation}
\label{meas1} y(t) = g(x(t),u(t)) + v(t)
\end{equation}
where $x$ is the weights vector in a neural network $g(x,u)$ that
maps the input vector $u$ to the output $y$. In (\ref{meas1}) we
denote by $u(t)$ the $t$'th input pattern and by $y(t)$ the corresponding
observed output pattern. The signal $v(t)$ represents either the
observation error or the modeling error.

Our aim is to estimate the weights vector $x$ from the
measurements $y(t), t=1,2,...,N$, given an initial estimate $x_0$.
Although the weights vector $x$ is constant (but unknown),
following \cite{inproceedings_Hassibi} we model it by
\begin{equation}
\label{scalsys} x(t+1) = x(t) + w(t)
\end{equation}
where $w(t)$ is the uncertainty in $x$ assumed to be a normally distributed white noise sequence with intensity $Q$, namely 
\begin{equation}
\label{RR}
E\{w(t)w(\tau)\} = Q \delta(t,\tau)
\end{equation}
where $\delta$ is the Kronecker delta function, namely attains the value of $1$ for $t = \tau$ and zero elsewhere, and where $Q$ may serve as a tuning parameter chosen according the trade off between fast convergence and accuracy. Similarly, $v$ is modeled also as a white noise sequence independent of $w$, where $E\{v(t)v(\tau)\} = R \delta(t,\tau)$ with $R$ being another parameter for tuning. Using a sequential
estimation scheme, we may want to estimate $x(k)$ by $\hat x(k)$ on the basis of accumulated measurements $Y(t)=\{y(\tau),\tau \in \{0,1,2,...,t\}\}$,
so that eventually at $t=N$ for large enough $N$, we will get that $\hat x(t) \approx
x$. We define the error signal at step $t$ by
\begin{equation}
\label{err} e(t) = x(t) - \hat x(t)
\end{equation}
and also define the following cost function
\begin{equation}
\label{cost}
\bar J = E\{e(t)e^T(t)\}
\end{equation}
Namely, we would like to minimize the Mean Square Estimation error covariance at each moment of the estimation process, namely at each training instance. At this point, one may wonder, whether the assumption of normal distribution of $w$ leading to a normally distributed weights vector $x$ can be justified. In our example we will see that this indeed the case. This assumption is supported to some extent in \cite{bellido1993}.  

{\bf Remark :} In the sequel we freely switch between e.g. $x(t)$ and $x_t$ meaning the same, when it serves the simplicity of representations.

\section{Training Algorithm}
It is well known (e.g. \cite{jazwinski2007stochastic}) that, 
\begin{equation}
\label{condme}
\hat x(t) = E(x|Y) :=\int_{-\infty}^{\infty} p(x|Y)dx
\end{equation}
where $p(x|Y)$ is the conditional probability density of $x$ given $Y$ and $E(x|Y)$ is the corresponding conditional expectation of the weights (i.e. state-vector) on the training examples outputs given by (1). Note that $\hat x$ is a function of $Y$. For the sake of completeness we further note that this result readily follows considering a general estimator of the form $\hat x = f(y)$. We would like to minimize the following loss function
$$
E\{ (x-f(y))^2\}=\int_{-\infty}^{\infty}p(y)\int_{-\infty}^{\infty} (x-f(y))^2 p(x|y)dx dy
$$
The result of (\ref{condme}) readily follows then from minimization of $\int_{-\infty}^{\infty} (x-f(y))^2 p(x|y)dx$.

Since (\ref{meas1}) defines a nonlinear relation between the weights and the training examples, we resort to using the Particle filtering algorithm that can deal with such a nonlinear situation without needing to calculate any gradients. As with the Kalman filter, the Particle filter \cite{Arulampalam02atutorial} employs two basic steps, the first one being a prediction step which takes the previous posterior (i.e. the conditional estimate $p(x_{t-1}|Y_{t-1})$ and using the state-vector evolution model of (\ref{scalsys}), namely, $p(x_t|x_{t-1})$, computes
\begin{equation}
    \label{predict}
    p(x_t|Y_{t-1})=\int p(x_t|x_{t-1}) p(x_{t-1}|Y_{t-1}) dx_{t-1}
\end{equation}
The second step, again similarly to the Kalman filter, is the measurement update step which is based on the Chapman - Kolmogorov identity
\begin{equation}
    \label{update}
    p(x_t|Y_{t})= P(x_t|Y_{t-1}) p(y_t|x_t)/p(y_t|Y_{t-1})
\end{equation}
where the denominator of the latter equation, is computed applying
\begin{equation}
\label{chapman}
p(y_t|Y_{t-1})=\int p(y_t|x_t) p(x_t | Y_{t-1})dx_t
\end{equation}
Equations (\ref{predict})-(\ref{chapman}) summarize the basic concepts of the Particle Filtering algorithm, where the probability density functions are represented by a finite number of samples of the state vector(also known as 'particles'), rather than the functions themselves. However, to make it work, some additional considerations such as importance sampling and re-sampling are to be employed. We note that we did not try to implement the algorithm of Particle Filtering, but rather used the one implemented within the control toolbox of $\mbox{MATLAB}^{TM}$. 

\section{Numerical Example}

We next apply Particle Filtering to the training of a functional neural network (FNN) (see e.g.\cite{noe2003} )which mimics the dynamics of the Henon Oscillator (see \cite{henon1976two}), out of noisy measurements of its outputs. The Henon Oscillator is described by the following nonlinear recursion :
\begin{equation}
\label{henon1} \xi(t)=c-a\xi(t-1)^2 +b \xi(t-2)
\end{equation}
where $a = 1.4,b=0.3$ and $c=1$. If one wants to identify the
system parameters $a,b$ and $c$ out of noisy measurements
\begin{equation}
\label{meas2} y(t) = \xi(t) + v(t) , k=0,1,2,...,N-1
\end{equation}
then the oscillator's model may be characterized as belonging to
the following class of models, where $C_{ij}$ are the unknown parameters which are to be found using a training algorithm :
\begin{eqnarray}
\label{henon2}
\xi(t) = c_{11} + c_{12} \xi(t-2) + c_{13} \xi(t-2)^2 \nonumber \\
+
         c_{21} + c_{22} \xi(t-1) + c_{23} \xi(t-1)^2 \nonumber
\end{eqnarray}
Defining the weights vector by
\begin{equation}
\label{henon3} x = col\{c_{11}+c_{21},c_{12},
c_{13},c_{22},c_{23}\}
\end{equation}
and the input vector by
$$
u(t) = \left[ \begin{array}{cc} \xi(t-1) & \xi(t-2) \end{array}
\right]^T
$$
we note that the output of the oscillator is given by (\ref{meas1}),
where
\begin{eqnarray}
\label{henon6} g(x,u(t)) = x_1 + x_2 u_2(t) + x_3 u_2(t)^2
\nonumber
\\ + x_4 u_1(t) + x_5 u_1(t)^2 :=C(u(t))x(t)
\end{eqnarray}
where
\begin{equation}
\label{henon6c} C(u(t)) = \left[ \begin{array}{ccccc} 1 & u_2(t) &
u_2(t)^2 & u_1(t) & u_1(t)^2 \end{array} \right]
\end{equation}
Noting, however, that having the noise corrupted outputs of the oscillator only, rather than its inner states,for training purposes, actually 
\begin{equation}
\label{meas4}
u(t) = \left[ \begin{array}{cc} y(t-1) & y(t-2) \end{array} \right]
\end{equation}
The training experiment is, therefore, based on running the Henon oscillator as a 'black-box' while we log its outputs $y(t)$ of (\ref{meas2}). The delayed versions of $y(t)$ are then used in (\ref{meas4} ) to form the network input vectors. Finally, the pairs $\{y(t),u(t)\}$ are the training examples that we use as inputs to the Particle Filter based training algorithm. We note that the same example was used in \cite{yaesh2014robust} using both a Kalman filter and $H_\infty$ filter. Here we use a Particle filter which is capable of training also general ANN rather than FNN, but we choose to test PF using this simple example, to verify that PF is not inferior to the Kalman filter and $H_\infty$-filter which are suited to deal with the linear case, and are not approximate solutions,in contrast to the PF which is not exact due to its finite number of particles.

The training results which were obtained using $Q = 0.5 \sqrt \Delta T$, where $\Delta T = 0.016$ and $R = 0.2$. The training results are depicted in Fig. 1, where we compare the time-history trained weights to the nominal true values. A satisfactory convergence of the weights to the true values is observed, following an initial transient. We next perform a simulation of the equation 
\begin{eqnarray}
\label{henon7}
\xi(t) = \hat x_1 + \hat x_2 \xi(t-2) + \hat x_3 \xi(t-2)^2 \nonumber \\
+
          \hat x_4 \xi(t-1) + \hat x_5 \xi(t-1)^2 \nonumber
\end{eqnarray}
The results of simulating the latter equation is shown in Fig. 2. clearly showing the classical Henon map.

\section{Conclusions}

A particle filtering based training concept for training general artificial neural networks, alleviating the difficulty involved in gradient based training processes such as back propagation is suggested. While testing the suggested approach to a deep learning example is necessary to establish its tractability, the suggested approach was demonstrated by considering a simple problem which has the same structure as general problems have, but due to its affinity in the weights, is solvable also using a Kalman filter and thus has a known global minimum. Therefore,the fact that the particle filter approach, in spite of being an approximate one (due to the finite number of particles), performed as good as the Kalman filter, should serve as an initial proof of concept and may encourage testing more involved training problems to the extent of dealing with deep learning problems. The fact that the particle filter may be implemented in a scalable manner \cite{part3} further motivates study of deep learning applications.

\bibliographystyle{unsrt}
\bibliography{particle1}

\begin{thebibliography}{10}

\bibitem{Haykin:1998:NNC:521706}
Simon Haykin.
\newblock {\em Neural Networks: A Comprehensive Foundation}.
\newblock Prentice Hall PTR, Upper Saddle River, NJ, USA, 2nd edition, 1998.

\bibitem{admm}
Taylor Gavin, Ryan Buremeister, Zheng Xu, Bharath Sing, Ankit Patel, and
  Goldstein Tom.
\newblock Training neural networks without gradients: A scalable admm approach.
\newblock In {\em arXiv:1605.02026v1 [cs.LG]}, 2016.

\bibitem{admm1}
Junxiang Wang, Fuxun Yu, Xiang Chen, and Liang Zhao.
\newblock Admm for efficient deep learning with global convergence.
\newblock In {\em KDD '19: Proceedings of the 25th ACM SIGKDD International
  Conference on Knowledge Discovery \& Data Mining.}, pages 39--44, 2019.

\bibitem{inproceedings_Hassibi}
Babak Hassibi, Ali H.~Sayed, and Thomas Kailath.
\newblock Lms is $h_\infty$ optimal.
\newblock In {\em ,}, volume~1, pages 74 -- 79 vol.1, 01 1994.

\bibitem{hassibi1995h}
Babak Hassibi and Thomas Kailath.
\newblock $h_\infty$ optimal training algorithms and their relation to
  backpropagation.
\newblock In {\em Advances in Neural Information Processing Systems}, pages
  191--198, 1995.

\bibitem{castillo1998introduction}
Enrique Castillo, Angel Cobo, JM~Guti{\'e}rrez, and E~Pruneda.
\newblock An introduction to functional networks with applications.
\newblock In {\em A Neural Based Paradigm}. Kluwer Academic Publishers, 1998.

\bibitem{yaesh2014robust}
Isaac Yaesh, Noelia~S{\'a}nchez Maro{\~n}o, and Uri Shaked.
\newblock Robust terminal filtering in an $h_\infty$ setting-application to
  affine functional networks training.
\newblock In {\em 2014 IEEE 28th Convention of Electrical \& Electronics
  Engineers in Israel (IEEEI)}, pages 1--5. IEEE, 2014.

\bibitem{jazwinski2007stochastic}
Andrew~H Jazwinski.
\newblock {\em Stochastic processes and filtering theory}.
\newblock Courier Corporation, 2007.

\bibitem{Arulampalam02atutorial}
M.~Sanjeev Arulampalam, Simon Maskell, and Neil Gordon.
\newblock A tutorial on particle filters for online nonlinear/non-gaussian
  bayesian tracking.
\newblock {\em IEEE TRANSACTIONS ON SIGNAL PROCESSING}, 50:174--188, 2002.

\bibitem{bellido1993}
I~Bellido and Fiesler E.
\newblock Do back propagation.
\newblock In {\em Proceedings of ICANN}, 1993.

\bibitem{noe2003}
Noelia S{\'a}nchez~Maro{\~n}o, Oscar Fontenla-Romero, Amparo Alonso-Betanzos,
  and Bertha Guijarro-Berdiñas.
\newblock Self-organizing maps and functional networks for local dynamic
  modeling.
\newblock In {\em ,}, pages 39--44, 01 2003.

\bibitem{henon1976two}
Michel H{\'e}non.
\newblock A two-dimensional mapping with a strange attractor.
\newblock In {\em The Theory of Chaotic Attractors}, pages 94--102. Springer,
  1976.

\bibitem{part3}
Pinalkumar Engineer, Rajbabu Velmurugan, and Sachin Patkar.
\newblock Scalable implementation of particle filter-based visual object
  tracking on network-on-chip (noc).
\newblock In {\em J Real-Time Image Proc 17, 1117–1134}, pages 1117--1134,
  2020.

\end{thebibliography}

\setlength{\baselineskip}{12pt}
\parskip 0.03cm

\newpage
\begin{figure}[h]
\centering
\includegraphics[width=\linewidth]{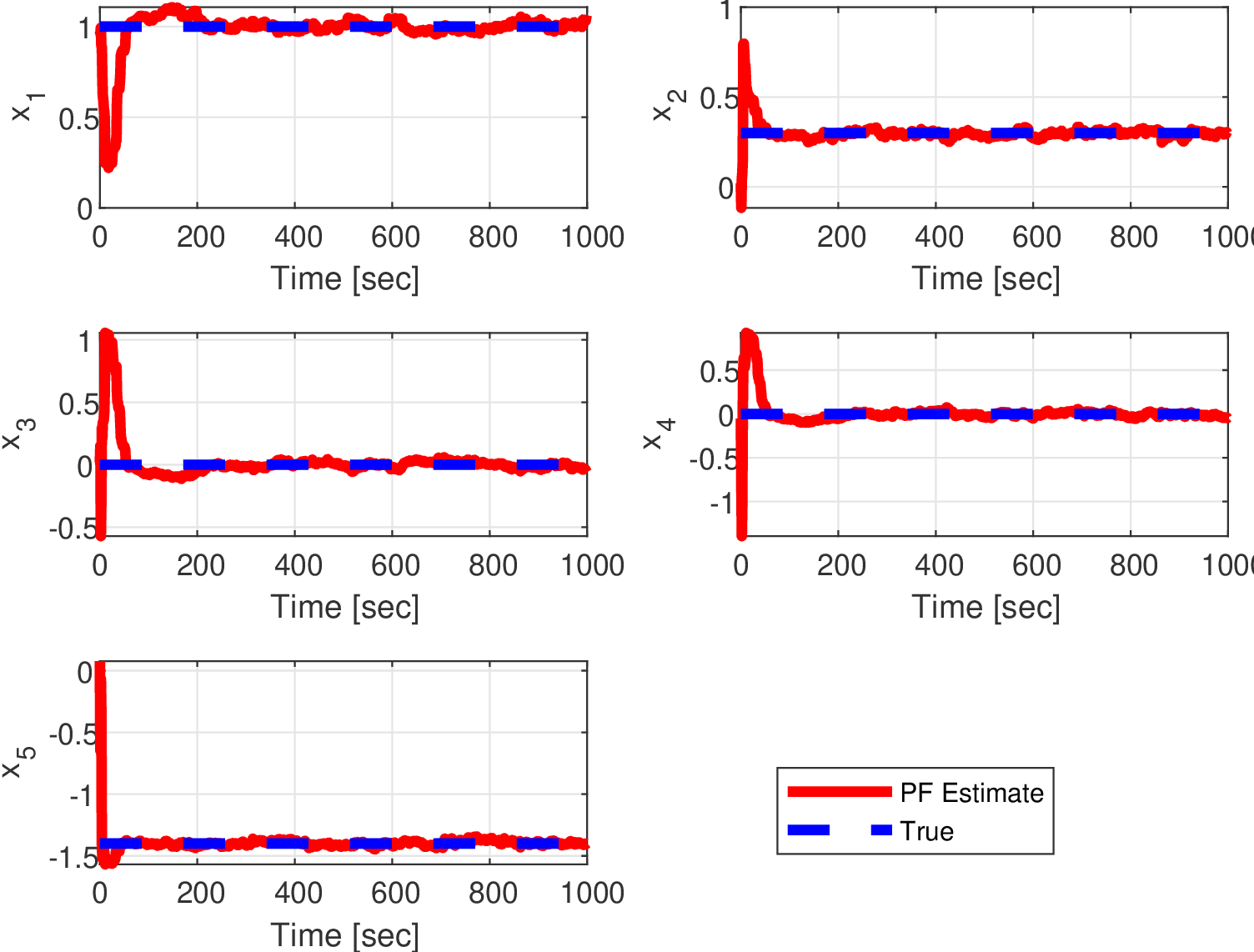}
\caption{Convergence of weights}
 \end{figure}

\begin{figure}[h]
\centering
\includegraphics[width=\linewidth]{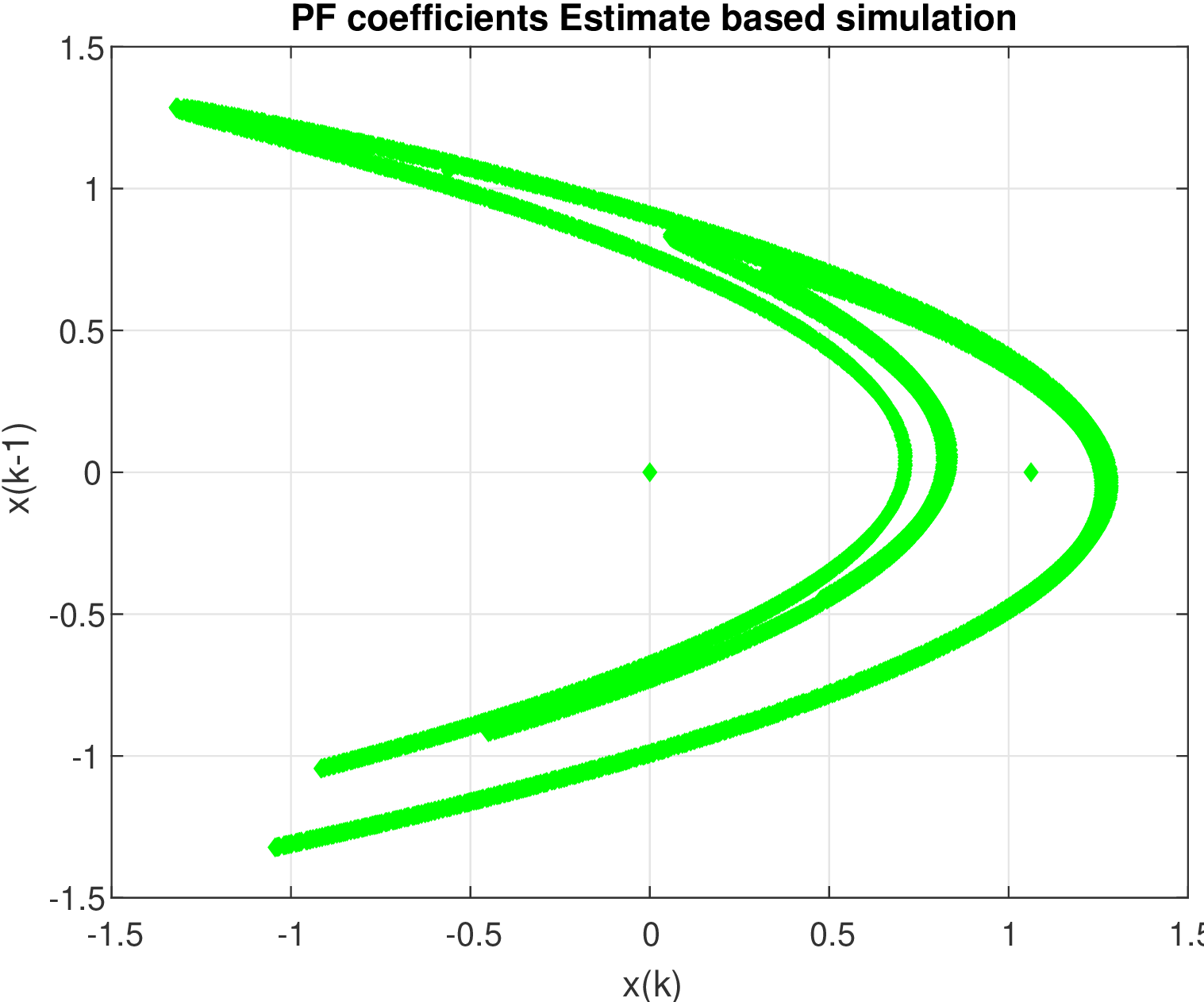}
\caption{Simulation results with training results for weights}
\end{figure}

\end{document}